\documentclass [12pt,a4paper,reqno]{amsart}
\usepackage{amsmath,amsthm,amscd}
\usepackage{amsfonts,amssymb}
\allowdisplaybreaks

\newcommand{\be}{\begin{equation}}
\newcommand{\ef}{\end{equation}}
\chardef\bslash=`\\ % p. 424, TeXbook
\newtheorem{thm}{Theorem}[section]
\newtheorem*{thm*}{Theorem}
\newtheorem{cor}[thm]{Corollary}
\newtheorem{lem}[thm]{Lemma}
\newtheorem{prop}[thm]{Proposition}

\theoremstyle{definition}

\newtheorem*{remark*}{Remarks}
\newtheorem*{defn*}{Definition}
\theoremstyle{remark}
\theoremstyle{remark*}

\numberwithin{equation}{section}

\newcommand{\G}{\Gamma}
\newcommand{\wt}{\widetilde}
\newcommand{\wh}{\widehat}
\newcommand{\bk}{\bigskip}
 \renewcommand{\sectionmark}[1]{}

\renewcommand{\Re}{\operatorname{Re}}
\newcommand{\iy}{\infty}
\newcommand{\Be}{Beltrami}
\newcommand{\hol} {holomorphic}
\newcommand{\qc} {quasiconformal}

\newcommand{\sh} {subharmonic}
\newcommand{\psh} {plurisubharmonic}
\newcommand{\ve}{\varepsilon}

\newcommand{\fc}{\frac}
\newcommand{\Te} {Teichm\"{u}ller}
\newcommand{\const}{\operatorname{const}}

 \usepackage{amsfonts}
\newcommand{\field}[1]{\mathbb{#1}}
\newcommand{\g}{\gamma}

\newcommand{\dl}{\delta}
\newcommand{\D}{\field{D}}
\newcommand{\om}{\omega}
\newcommand{\z}{\zeta}
\newcommand{\ov}{\overline}
\newcommand{\vp}{\varphi}
\newcommand{\hC}{\widehat{\field{C}}}
\newcommand{\C}{\field{C}}

\newcommand{\B}{\mathbf{B}}
\newcommand{\T}{\mathbf{T}}

\newcommand{\Hol}{\operatorname{Hol}}
\newcommand{\id}{\operatorname{id}}
\newcommand{\grad}{\operatorname{grad}}

\newcommand{\uTs} {universal Teichm\"{u}ller space}
\newcommand{\Belt} {\operatorname{Belt}}
\newcommand{\Ca} {Carath\'{e}odory}
\newcommand{\Gr} {Grunsky}
\newcommand{\Ko} {Kobayashi}
\newcommand{\Om} {\Omega}
\newcommand{\vk} {\varkappa}
\newcommand{\kp} {\kappa}
\newcommand{\x} {\mathbf x}
\renewcommand{\a} {\alpha}
\newcommand{\ld}{\lambda}
\newcommand{\Ld}{\Lambda}

\begin{document}

\title{Complex geodesics and variational calculus for univalent
functions}

\author{Samuel L. Krushkal}

\begin{abstract}It turns out that complex geodesics in Teichm\"{u}ller
spaces with respect to their invariant metrics are intrinsically
connected with variational calculus for univalent functions.

We describe this connection and show how geometric features associated
to these metrics and geodesics provide deep distortion results for
various classes of functions with quasiconformal extensions and
create new phenomena which do not appear in the classical geometric
function theory.
\end{abstract}

\date{\today\hskip4mm ({\tt geodvarCADS7rev.tex})}

\maketitle

\bigskip
{\small {\textbf {2010 Mathematics Subject Classification:} Primary:
30C55, 30C62, 30C75, 30F60, 32F45; Secondary: 30F45, 46G20}

\medskip


\markboth{Samuel Krushkal}{Complex geometry and variational
calculus} \pagestyle{headings}

\bigskip\bigskip
This paper is an extended version of my talk at the conference
``Complex Analysis and Dynamical Systems VII". Our main goal is to
reveal the intrinsic connection between the complex geodesics in \Te
\ spaces and extremals of the basic functionals in geometric
function theory. This connection causes new phenomena that do not
arise in the classical theory and provides an alternative approach
to variational problems as well.

\bk
\section{Key theorems on invariant metrics and geodesics
on \Te \ spaces}

\bk\subsection{\Te \ spaces}

We shall use the notations $\hC = \C \cup \{\iy\}$,
$$
\D = \{|z| < 1\}, \ \D^* = \{z \in \hC: \ |z| > 1\}
$$
and consider two \Te \ spaces: the \uTs \ $\T = \T(\D)$ and the
 \Te \ space $\T_1 = \T(\D^0)$ of the punctured disk
$\D^0 = \D \setminus \{0\}$.

Every \Te \ space $\wt \T$ is a complex Banach manifold, thus it
possesses the invariant {\bf Cara\-th\'{e}odory} and {\bf \Ko} \
distances (the
smallest and the largest among all holomorphically non-expanding
metrics).

Denote these metrics by $c_{\wt \T}$ and $d_{\wt \T}$, and
let $\tau_{\wt \T}$ be the intrinsic \Te \ metric of this space
canonically determined by \qc \ maps. Their
infinitesimal Finsler forms are defined on the tangent bundle
$\mathcal T \wt \T$ of $\wt \T$), and satisfy
$$
c_{\wt \T}(\cdot, \cdot) \le d_{\wt \T}(\cdot, \cdot) \le \tau_{\wt
\T}(\cdot, \cdot).
$$
By the fundamental Royden-Gardiner theorem, the metrics $d_{\wt \T}$
and $\tau_{\wt \T}$ (and their infinitesimal forms) are equal; see,
e.g. \cite{EKK}, \cite{EM}, \cite{GL}, \cite{Ro}.

\bk So an important and still open question is how about the \Ca \
metric of $\wt T$: does it also coincide with the Teichm\"{u}ller
metric? Our first theorem yields a positive answer for both spaces
$\T$ and $\T_1$. Let $\wt T$ denote either from these spaces.

\bk
\begin{thm} The \Ca \ metric of the space $\wt \T$ coincides with
its \Ko \ metric, hence all invariant non-expanding metrics on
$\wt \T$ are equal its \Te \ metric, and
$$
 c_{\wt T}(X_1, X_2) = d_{\wt T}(X_1, X_2) = \tau_{\wt T}(X_1, X_2)
= \inf \{d_\D(h^{-1}(\vp), h^{-1}(\psi)): \ h \in \Hol(\D, \wt \T)\},
$$
where $d_\D$ denotes the hyperbolic metric of the unit disk of
curvature $-4$.

Similarly, the infinitesimal forms of these metrics coincide with
the Finsler metric $F_{\wt \T}(\vp, v)$ generating $\tau_{\wt \T}$
and have \hol \ sectional curvature $-4$.
\end{thm}

\bk This was known only for the \uTs \ and underlies various
applications. A new proof for this space is given in \cite{Kr8}. Its
arguments (based on \psh \ features and geometry of the \Gr \
operator on univalent functions) can be extended to the second space
$\T_1$. Since Theorem 1.1 underlies all other results of this paper,
we present here its proof in the lines of \cite{Kr8} focusing mainly
on the new case $\wt \T = \T_1$.

Recall that the \uTs \ $\T =\T(\D)$ is the space of quasisymmetric
homeomorphisms of the unit circle $S^1 = \partial \D$ quotiented
by M\"{o}bius maps.
The canonical complex Banach structure on $\T$ is defined by the
quotient of the ball of \Be \ coefficients
$$
\Belt(\D)_1 = \{\mu \in L_\iy(\C) : \ \mu|\D^* = 0, \ \| \mu \| <
1\},
$$
letting $\mu, \nu \in \Belt(\D)_1$ be equivalent if the
corresponding \qc \ maps $f^\mu, f^\nu \in \Sigma^0$
coincide on $S^1$ and passing to Schwarzian derivatives
$$
S_w(z) = \Bigl(\fc{w^{\prime\prime}(z)}{w^\prime(z)}\Bigr)^\prime
- \fc{1}{2} \Bigl(\fc{w^{\prime\prime}(z)}{w^\prime(z)}\Bigr)^2
\quad (z \in \D^*, \ \ w = f^\mu|\D^*).
$$
Here $\Sigma^0$ denotes the collection of univalent functions
$$
f(z) = z + b_0 + b_1 z^{-1} + \dots
$$
in $\D^*$ admitting \qc \ extensions to the unit disk $\D$ with
$f(1) = 1$; \ $f^\mu$ denotes the solution of the \Be \ equation
$\ov{\partial} w = \mu \partial f$ on $\C$ with $\mu \in \Belt(\D)_1$
and above normalization.

The Schwarzians $S_{f^\mu}$ run over a bounded domain in the
Banach space $\B$ of hyperbolically bounded \hol \ functions on $\D^*$
with norm
$$
\|\vp\| = \sup_{\D^*} (|z|^2 - 1)^2 |\vp(z)|.
$$
This domain models the space $\T$, and the defining projection
$\phi_\T: \mu \to S_{f^\mu}$ is \hol \ as a map $L_\iy(\D) \to \B$.

The intrinsic \textbf{\Te \ metric} of the space $\T$ is given by
$$
\tau_\T (\phi_\T (\mu), \phi_\T (\nu)) = \frac{1}{2} \inf \bigl\{
\log K \bigl( w^{\mu_*} \circ \bigl(w^{\nu_*} \bigr)^{-1} \bigr) : \
\mu_* \in \phi_\T(\mu), \nu_* \in \phi_\T(\nu) \bigr\},
$$
where $\phi_\T$ is the factorizing \hol \ projection $\Belt(\D)_1 \to
\T$ and
$$
K(w^\mu) = (1 + \|\mu\|_\iy)/(1 - \|\mu\|_\iy)
$$
denotes the maximal \qc \ dilatation of the map $w^\mu$.
This metric is the integral form of the infinitesimal Finsler
metric
$$
  F_\T(\phi_\T(\mu), \phi_\T^\prime(\mu) \nu) = \inf
\{\|\nu_*/(1 - |\mu|^2)\|_\iy: \ \phi_\T^\prime(\mu) \nu_* =
\phi_\T^\prime(\mu) \nu\}
$$
on the tangent bundle $\mathcal T\T$ of $\T$.

The space $\T_1 =\T(\D^0)$ is canonically isomorphic to the subspace
$\T(\G) = \T \cap \B(\G)$, where $\G$ is a cyclic parabolic Fuchsian
group in $\D$ uniformizing the punctured disk $\D^0$ and
$$
\B(\G) = \{\vp \in \B: \ (\vp \circ \g) (\g^\prime)^2 = \vp, \ \
\g \in \G\} \quad \text{in} \ \ \D^*.
$$
By \cite{Be}, this space is biholomorphically equivalent to the Bers
fiber space $\mathbf F(\T)$ over \uTs, which underlies the proof of
Theorem 1.1 for $\T_1$.

\bk
\subsection{Complex geodesics}
If $X$ is a domain in a complex Banach space $E$ endowed
with a pseudo-distance $\rho_X$, then a \hol \ map $h: \ \D \to X$
is called a {\bf complex $\rho$-geodesic} if there exist $t_1 \ne
t_2$ in $\D$ such that
$$
d_\D(t_1, t_2) = \rho_X(h(t_1), h(t_2));
$$
one says also that the points $h(t_1)$ and $h(t_2)$ can be joined by
a complex $\rho$-geodesic (see \cite{Ve}).

If $h$ is a complex $c_X$-geodesic then it also is $d_X$-geodesic
and the above equality holds for all points $t_1, \ t_2 \in \D$, so
$h(\D)$ is a \hol \ disk in $X$ hyperbolically isometric to $\D$.

\bk The second basic theorem follows from Theorem 1.1 and the
properties of extremal \qc \ maps.

\bk
\begin{thm} (i) Any two points of the space $\wt \T$ can be joined by
a complex geodesic. The geodesic joining a Strebel point with the
base point is unique and defines the corresponding \Te \ extremal
disk.

(ii) For any point $\vp \in \wt \T$ and any nonzero tangent vector
$v$ at this point, there exists at least one complex geodesic $h: \
\D \to \wt \T$ such that $h(0) = \vp$ and $h^\prime(0)$ is collinear
to $v$.
\end{thm}

\subsection{}
The coincidence of invariant metrics and existence of complex
geodesics was known only for convex domains in dual Banach spaces
(see \cite{DTV}). Duality is needed to apply the Alaoglu-Banach
theorem.

Theorem 1.1 implies, together with the definition of complex
geodesics, that these geodesics in $\wt \T$ are the \Te \ geodesic
disks in this space.

Such disks are uniquely determined for Strebel points; on the other
hand, Tanigawa constructed in \cite{Ta} the extremal
\Be \ coefficients $\mu_0$ with nonconstant $|\mu_0(z)| <
\|\mu_0\|_\iy$ on a set of positive measure for which there exist
infinitely many distinct geodesic segments in the \uTs \ $\T$
joining the points $\phi_\T(\mathbf 0)$ and $\phi_\T(\mathbf
\mu_0)$. All these segments belong to different complex geodesics
joining the indicated points.

\bk
\section{Proof of Theorem 1.1}

The proof of this underlying theorem involves essentially
the properties of the Grunsky coefficients of univalent functions.
As an important consequence of this theorem, one obtain Theorem 1.2
on complex geodesics of these metrics which will be applied
to variational calculus.

\bk\noindent
$\mathbf{1^0}$. {\bf Auxiliary lemmas}.
Consider the space $\Sigma^0$ of the
univalent $\hC$-\hol \ functions $f(z) = z + b_0 + b_1 z^{-1} +
\dots$ on the disk $\D^*$ with \qc \ extension to $\hC$ satisfying
$f(1) = 1$. Their \Gr \ coefficients $\a_{mn}$ are defined via
 \be\label{2.1}
\log \frac{f(z) - f(\z)}{z - \z} = - \sum\limits_{m,n = 1}^{\iy}
\a_{mn} z^{-m} \z^{-n}, \quad (z, \z) \in (\D^*)^2
\end{equation}
choosing the branch of the logarithmic function which vanishes as $z
= \z \to \iy$. The quantity
$$
\vk(f) := \sup \Big\{ \Big\vert \sum\limits_{m,n=1}^\iy \sqrt{mn} \
\a_{mn} x_m x_n \Big\vert \Big\}
$$
where the supremum is taken over the points of the unit sphere
$S(l^2) = \{\|\mathbf x \| = 1\}$ in the Hilbert space of sequences
$\mathbf x = (x_n)$ with $\|\mathbf x \|^2 = \sum\limits_{1}^{\iy}
|x_n|^2$, is called the {\Gr \ norm} of $f$, and by a result of
Grunsky \cite{Gr} the inequality $\vk(f) \le 1$ is the necessary and
sufficient condition for univalence of $f$ in the disk $\D^*$.

The coefficients $\a_{mn}$ are polynomials of the initial Taylor
coefficients $b_1, b_2, \dots, b_{m+n_1}, \dots$ of $f$, hence
depend holomorphically on \Be \ differentials $\mu \in \Belt(\D)_1$
and on the Schwarzian derivatives $S_f \in \B$ of these functions
(see, e.g., \cite{GL}, \cite{Kr1}). This generates for a fixed $\x =
(x_n) \in S(l^2)$ the \hol \ map
 \be\label{2.2}
 h_\x(S_f) = \sum\limits_{m,n=1}^\iy \ \sqrt{m n}
\ \a_{m n} (S_f) \ x_m x_n : \ \T \to \D,
\end{equation}
and $\sup_{\x\in S(l^2)} |h_\x(S_f)| = \vk(f)$. Both the \Te \ and
\Gr \
norms of $f$ are continuous \psh \ functions of $S_f$ on $\T$ (see,
e.g. \cite{Kr4}, \cite{Kr7}).

Note that the convergence and holomorphy of the series (2.1) simply
follow from the inequalities
$$
\Big\vert \sum\limits_{m=j}^M \sum\limits_{n=l}^N \ \sqrt{m n} \
\a_{mn} x_m x_n \Big\vert^2 \le \sum\limits_{m=j}^M |x_m|^2
\sum\limits_{n=l}^N |x_n|^2
$$
(for any finite $M, \ N$) which, in turn, are a consequence of the
classical area theorem (see, e.g., [30, p. 61]).

\bigskip
Let $A_1(\D)$ be the subspace of $L_1(\D)$ formed by \hol \
functions in $\D$, and
$$
A_1^2 = \{\psi \in A_1(\D): \ \psi = \om^2\};
$$
this set consists of the integrable \hol \ functions \ on $\D$
having only zeros of even order. In fact, all $\psi \in A_1^2$ are
of the form
$$
\psi(z) = \fc{1}{\pi} \sum\limits_{m+n=2}^\iy \ \sqrt{m n} \ x_m x_n
z^{m+n-2},
$$
with $\|\x\|_{l^2} = \|\om\|_{L_2}$. Put
$$
<\mu, \psi>_\D = \iint_\D  \mu(z) \psi(z) dx dy, \quad \mu \in
L_\iy(\D), \ \psi \in L_1(\D) \ \ (z = x + iy)
$$
and
$$
\a_\D(f) = \sup \ \{|\langle \mu_0, \psi\rangle_\D|: \ \psi \in
A_1^2, \ \|\psi\|_{A_1(\D)} = 1\}.
$$

The \Gr \ norm of every $f \in \Sigma^0$ is dominated by its \Te \
norm, i.e., $\vk(f) \le k(f)$ \cite{Ku1}. We shall need a stronger
estimate given by

\begin{lem} \cite{Kr2}, \cite{Kr7} For all $f \in \Sigma^0$,
$$
\vk(f) \le k \fc{k + \a_\D(f)}{1 + \a_\D(f) k}, \quad k = k(f),
$$
and $\vk(f) < k$ unless
 \be\label{2.3}
\a_\D(f) = \|\mu_0\|_\iy,
\end{equation}
and the last equality is equivalent to $\vk(f) = k(f)$. Moreover,
for small $\|\mu\|_\iy$,
$$
\vk_r(f) = \sup \ |\langle \mu, \psi\rangle_\D| + O(\|\mu\|_\iy^2),
\quad \|\mu\|_\iy \to 0,
$$
with the same supremum as in (2.3).

If $\vk(f) = k(f)$ and the equivalence class of $f$ (the collection
of maps equal $f$ on $S^1 = \partial \D^*$) is a Strebel point, then
$\mu_0$ is necessarily of the form
 \be\label{2.4}
\mu_0 = \|\mu_0\|_\iy |\psi_0|/\psi_0 \ \ \text{with} \ \ \psi_0 \in
A_1^2.
\end{equation}
\end{lem}

Geometrically, the equality (2.3) means the equality of the \Ca \
and \Te \ distances on the geodesic disk
$$
\D(\mu_0) = \{\phi_\T(t\mu_0 /\|\mu_0\|): t \in \D\} \subset \T
$$
(compare with Kra's theorem for finite dimensional \Te \ spaces
\cite{K}).

For functions $f \in \Sigma^0$ \hol \ in the closed disk
$\ov{\D^*}$, the relation (2.4) was also obtained (by a different
method) in \cite{Ku3}.

\bigskip
One defines for $f \in \Sigma^0$ the complex homotopy
$$
f_t (z) = t f(z/t) =  z + b_0 t + b_1 t^2 z^{-1} + b_2 t^3 z^{-2} +
...: \ \D^* \times \D \to \hC.
$$
Then $S_{f_t}(z) = t^{-2} S_f(t^{-1} z)$, and moreover, this
point-wise map determines a \hol \ map
 \be\label{2.5}
\beta_f(t) =  S_{f_t}(\cdot): \ \D \to \B.
\end{equation}
The homotopy disks $\D(S_f) = \beta_f(\D)$ have only cuspidal
critical points (are branched simply connected Riemann surfaces
placed in $\T$) and foliate the space $\T$.

Each homotopy map $f_t$ admits an extremal extension to $\D$ of \Te
\ type with dilatation $k(f_t) \le k(f) |t|^2$. This bound is sharp
and is improved for
$$
f(z) = z + b_0 + b_m z^{-m} + \dots \quad (m > 1, \ b_m \ne 0),
$$
via $k(f_t) \le k(f) |t|^{m+1}$. Due to \cite{KK2},
$$
k(f_t) = \fc{m + 1}{2} |b_m| |t|^{m+1} + O(t^{m+2}), \quad t \to 0.
$$

The following lemma belongs to K\"{u}hnau \cite{Ku3}, see also [24,
Section 6].

\begin{lem} In the case $b_1 \ne 0$, there holds for sufficiently
small $|t| \le r_0(f)$ the equality $\vk(f_t) = k(f_t)$.
\end{lem}

Indeed, due to \cite{Ku3}, for small $|t|$ the extremal \qc \
extension of $f_t$ to $D$ is defined by a nonvanishing \hol \
quadratic differential. Generically, $r_0(f) < 1$, which is
connected with critical points of the homotopy disk.

Note that $S_f(z) = - 6 b_1 z^{-4} + O(z^{-5})$, so $b_1 = 0$ only
for $f$ with
 \be\label{2.6}
\lim\limits_{z\to \iy} z^4 S_f(z) = - 6 b_1 = 0.
\end{equation}
It suffices to prove the theorem for the set of $S_f \in \T$
such that $ b_1 \ne 0$ which is open and dense in $\T$.

The following lemma is a special case of a local existence theorem
from \cite{Kr1}.

\bk
\begin{lem}Let $D$ be a simply connected domain on the Riemann sphere
$\hC$. Assume that there are a set $E$ of positive two-dimensional
Lebesgue measure and a finite number of points
 $z_1, z_2, ..., z_m$ distinguished in $D$. Let
$\a_1, \a_2, ..., \a_m$ be non-negative integers assigned to $z_1,
z_2, ..., z_m$, respectively, so that $\a_j = 0$ if $z_j \in E$.

Then, for a sufficiently small $\ve_0 > 0$ and $\varepsilon \in (0,
\varepsilon_0)$, and for any given collection of numbers $w_{sj}, s
= 0, 1, ..., \a_j, \ j = 1,2, ..., m$ which satisfy the conditions
$w_{0j} \in D$, \
$$
|w_{0j} - z_j| \le \ve, \ \ |w_{1j} - 1| \le \ve, \ \ |w_{sj}| \le
\ve \ (s = 0, 1, \dots   a_j, \ j = 1, ..., m),
$$
there exists a \qc \ self-map $h$ of $D$ which is conformal on $D
\setminus E$ and satisfies
$$
h^{(s)}(z_j) = w_{sj} \quad \text{for all} \ s =0, 1, ..., \a_j, \ j
= 1, ..., m.
$$
Moreover, the \Be \ coefficient $\mu_h(z) = \partial_{\bar z}
h/\partial_z h$ of $h$ on $E$ satisfies $\| \mu_h \|_\iy \leq M
\ve$. The constants $\ve_0$ and $M$ depend only upon the sets $D, E$
and the vectors $(z_1, ..., z_m)$ and $(\a_1, ..., \a_m)$.

If the boundary $\partial D$ is Jordan or is $C^{l + \a}$-smooth,
where $0 < \a < 1$ and $l \geq 1$, we can also take $z_j \in
\partial D$ with $\a_j = 0$ or $\a_j \leq l$, respectively.
\end{lem}

The following lemma relies on classical results on compactness in
the dual week$^*$ topology and is a consequence of the
Alaoglu-Bourbaki theorem.

\begin{lem} Let $Y$ be reflexive and let $\Om$ be a bounded set
in  $\Hol(G, Y)$. Then any sequence from $\Om$ contains a
subsequence which weakly converges to a \hol \ map from $G$ to $Y$.
\end{lem}

Hence, the compactness (strong or weak) is actually required only
for the image of $G$ in $Y$ under a given family of maps.

In our case $Y = \C$, and for each $x \in \wt \T$ its orbit $h(x)$
is located in the unit disk which is compact. This implies that for
any point $x_0$ from $\wt \T$ there exists a \hol \ map $h_0: \ \wt
\T \to \D$ with $h_0(0) = 0$ and $d_\D(0, h_0(x_0)) = c_\T(\mathbf
0, x_0)$.

\bk\noindent
$\mathbf{2^0}$. {\bf Case $\wt \T = \T$}. Consider for the
Schwarzians $S_f \in \T$ their homotopy disks $\Delta(S_f)$.
It follows from Lemmas 2.1 and 2.2 that for small $|t| \le r_0(f)$,
 \be\label{2.7}
c_\T(\mathbf 0, S_{f_t}) = d_\T(\mathbf 0, S_{f_t}) =
\tau_\T(\mathbf 0, S_{f_t}) = \tanh^{-1} |t|,
\end{equation}
and similarly for the corresponding infinitesimal metrics. Our goal
now is to extend this equality to all points of $\Delta(S_f)$.

We first illustrate the arguments on the \uTs \ $\T$. The proof
for the space $\T_1 = \T(\D^0)$ will be given separately in
the next step.

Consider in the tangent bundle $\mathcal T(\T) = \T \times \B$ the
\hol \ disks $\wh \D(S_f)$ over the homotopy disk of $f$ formed by
the points $(\vp_t, v)$, where $v = \phi_\T^\prime[\vp_t] \mu \in
\B$ is a tangent vector to $\T$ at the point $\vp_t = S_{f_t}$, and
$\mu$ runs over the ball
$$
\Belt(D_{\vp_t})_1 = \{\mu \in L_\iy(\C): \ \mu|D_{\vp_t}^* = 0, \
\|\mu\|_\iy < 1\}.
$$
Here $D_\vp$ and $D_\vp^*$ denote the images of $\D$ and $\D^*$
under $f = f_\vp \in \Sigma^0$ with $S_f = \vp$.

Pick the unit tangent vectors $v$ to $\T$ at the origin and \hol \
maps $\wh h_v: \ \T \to \D$ with $\wh h_v(\mathbf 0) = 0$ and
maximal lengths $|d \wh h_v(\mathbf 0) v| = \mathcal C_\T(\mathbf 0,
v)$ (given by Lemma 2.4) and apply these maps to pulling back the
hyperbolic metric on the unit disk
$$
ds = |dt|/(1 - |t|^2)
$$
of Gaussian curvature $- 4$ onto the disks $\wh \D(S_f)$. Then we
obtain on these disks the logarithmically \sh \ metrics $ds = \wh
\ld_{\wh h_v}(t) |dt|$ with
$$
\wh \ld_{\wh h_v}(t) = \wh h_v^* \ld_\D = \fc{|\wh h_v^\prime(t)|}{1
- |\wh h_v(t)|^2},
$$
having Gaussian curvature $- 4$ at noncritical points. The upper
envelope of these metrics
$$
\ld_0(t) =\sup_{\wh h_v} \wh \ld_{\wh h_v}(t),
$$
followed by its upper semicontinuous regularization depends only on
the points $\vp \in \T$, hence it descends to a logarithmically \sh
\ metric $\ld_0$ on the underlying disk $\D(S_f)$. The underlying
metric $\ld_0$ has at each of its noncritical point $t_0$ a supporting
metric of curvature $- 4$, thus the generalized Gaussian curvature
of $\ld_0$ satisfies $\kp_{\ld_0} \le - 4$, or equivalently,
 \be\label{2.8}
 \Delta \log \ld_0 \ge 4 \ld_0^2.
\end{equation}
Here $\Delta$ means the generalized Laplacian
$$
\Delta \ld(t) = 4 \liminf\limits_{r \to 0} \frac{1}{r^2} \Big\{
\frac{1}{2 \pi} \int_0^{2\pi} \ld(t + re^{i \theta}) d \theta -
\ld(t) \Big\} \quad (0 \le \ld(t) < \iy).
$$

Now the metrics $\ld_{\mathcal C}$ and $\ld_{\mathcal K}$ can be
compared on the disk $\D(S_f)$ by Minda's maximum principle for
solutions of (2.8) given by

\begin{lem} \cite{Mi} If a function $u : \ \Om \to [- \iy,
+ \iy)$ is upper semicontinuous in a domain $\Om \subset \C$ and its
generalized Laplacian satisfies the inequality $\Delta u(z) \ge K
u(z)$ with some positive constant $K$ at any point $z \in \Om$,
where $u(z) > - \iy$, and if $\limsup\limits_{z \to \z} u(z) \le 0$
for all $\z \in \partial \Om$, then either $u(z) < 0$ for all $z \in
\Om$ or else $u(z) = 0$ for all $z \in \Om$.
\end{lem}

It is applied to the ratio $u = \log(\ld_{\mathcal C}/\ld_{\mathcal
K}) = \log \ld_{\mathcal C} - \log \ld_{\mathcal K}$. One needs to
consider only the points $S_{f_t}$ with $|t| > r_0$. In view of the
upper semicontinuity of both metrics the set $U_0$ of points $t \in
\D, \ |t| > r_0$, where $\ld_{\mathcal C}(t) < \ld_{\mathcal K}(t)$,
is open. Arguing similar to \cite{Mi}, one derives from (2.6) and
Lemma 2.5 that $\ld_{\mathcal C} = \ld_{\mathcal K}$ on $U_0$.
Continuing in a similar way, one extends this equality first to all
noncritical points of the disk $\D(S_f)$ and then by continuity of
metrics to the whole homotopy disk. This also provides the equality
of the global distances $c_\T$ and $d_\T$ on this disk.

Now, by applying the infinitesimal analog of Montel's normality
theorem and Lemma 2.4, one gets, letting $t \to 1$,
 \be\label{2.9}
c_\T(\mathbf 0, S_f) = d_\T(\mathbf 0, S_f).
\end{equation}
The case of two arbitrary points $\vp_1, \vp_2 \in T$ is reduced to
(2.9) in a standard way using the right translations of
$\Belt(\D)_1$. The details will be given in the next step.

\bk\noindent
{\bf Remark}. We have distinguished in the above proof the canonical
disk $\D^*$ as the base point of $\T$. In fact, the arguments work
for any other base point (quasidisk $D^*$) of this space.
Accordingly, one can define by (2.1) the generalized Grunsky
coefficients $\a_{m n}(f)$ for univalent functions $f(z)$ in this
quasidisk $D^*$ and construct similar to (2.2) the holomorphic maps
$h_\x(S_f): \ \T \to \D$ moving the base point $D^*$ into the
origin (see \cite{Kr8}).

\bk\noindent
$\mathbf{3^0}$. {\bf Proof for the space $\T_1$}.
First recall that the elements of the space
$\T_1 = \T(\D \setminus \{0\})$
are the equivalence classes of \Be
\ coefficients $\mu \in \Belt(\D)_1$ so that the corresponding \qc \
automorphisms $w^\mu$ of the unit disk coincide on both boundary
components (unit circle $S^1$ and the puncture $z = 0$)
and are homotopic on $\D \setminus \{0\}$. This space can be endowed
with a canonical complex structure of a complex Banach manifold and
embedded into $\T$ using uniformization.

Namely, the punctured disk $\D \setminus \{0\})$ is conformally
equivalent to the factor
$\D/\G$, where $\G$ is a cyclic parabolic Fuchsian group acting
discontinuously on $\D$ and $\D^*$. The functions $\mu \in
L_\iy(\D)$ are lifted to $\D$ as the \Be \ $(-1, 1)$-measurable
forms  $\wt \mu d\ov{z}/dz$ in $\D$ with respect to $\G$, i.e., via
$(\wt \mu \circ \g) \ov{\g^\prime}/\g^\prime = \wt \mu, \ \g \in
\G$, forming the Banach space $L_\iy(\D, \G)$.

We extend these $\wt \mu$ by zero to $\D^*$ and consider the unit ball
$\Belt(\D, \G)_1$ of $L_\iy(\D, \G)$. Then the corresponding
Schwarzians $S_{w^{\wt \mu}|\D^*}$ belong to $\T$. Moreover, $\T_1$
is canonically isomorphic to the subspace $\T(\G) = \T \cap \B(\G)$,
where $\B(\G)$ consists of elements $\vp \in \B$ satisfying $(\vp
\circ \g) (\g^\prime)^2 = \vp$ in $\D^*$ for all $\g \in \G$.

Due to the Bers isomorphism theorem, the space $\T_1$ is
biholomorphically equivalent to the Bers fiber space
$$
\mathcal F(\T) = \{\phi_\T(\mu), z) \in \T \times \C: \ \mu \in
\Belt(\D)_1, \ z \in w^\mu(\D)\}
$$
over the \uTs \ with \hol \ projection $\pi(\psi, z) = \psi$ (see
\cite{Be}). This fiber space is a bounded domain in $\B \times \C$.

Note that now the admissible quadratic differentials defining
\Te \ extremal coefficients $\mu_0 \in \Belt(\D)_1$ must be integrable
and \hol \ only on $\D \setminus \{0\}$, thus can have
simple pole at $z = 0$, i.e., $\mu_0 = k |\psi_0|/\psi_0$ with
$$
\psi_0(z) = c_{-1} z^{-1} + c_0 + c_1 z + \dots \ , \quad 0 < |z| <
1.
$$

To prove the theorem for the space $\T_1$, we establish the
equality of invariant metrics on this fiber space.
We again model the space $\T$ as a domain in the space $\B$ formed
by the Schwarzians $S_{f^\mu}$ of functions $f^\mu \in \Sigma^0$.

We associate with $f^\mu$ the odd function
\be\label{2.10}
\mathcal R f^\mu(z) := (f^\mu(z^2) - f^\mu(0))^{1/2} = z +
\fc{ b_0 - f^\mu(0)}{2 z} + \fc {b_3^\prime}{z^3} + \dots
\end{equation}
whose \Gr \ coefficients $\a_{m n}(\mathcal R f^\mu)$ are
represented as polynomials of the first Taylor coefficients of the
original function $f^\mu$ and of $a = f^\mu(0)$. Hence,
$\a_{m n}(\mathcal R f^\mu)$ depend holomorphically on the
Schwarzians $\vp = S_{f^\mu} \in \T$ and on values $f^\mu(0)$, i.e.,
on pairs $X = (\vp, a)$ which are the points of the fiber space
$\mathcal F(\T)$. This joint holomorphy follows from Hartog's
theorem on separately \hol \ functions extended to Banach domains.

The square root transform act on quadratic differentials
$\psi_0 dz^2$ via
$$
\mathcal R_{*} \psi_0 = \psi_0(z^2) 4 z^2 dz^2 = (4 c_{-1} + \wt c_2
z^2 + \dots) dz^2;
$$
so one can apply to $\mathcal R f^\mu$ Lemmma 2.2.

We construct by (2.2) for $\mathcal R f^\mu$ the
corresponding \hol \ functions
$$
\wh h_\x(S_f) = \sum\limits_{m,n=1}^\iy \ \sqrt{m n}
\ \a_{m n} (S_{\mathcal R f}) \ x_m x_n : \ \T \to \D
$$
mapping the domain $\mathcal F(\T)$ to the unit disk and
satisfying
 \be\label{2.11}
\sup_{\x\in S(l^2)} \wh h_\x(S_f) = \vk(\mathcal R f).
\end{equation}

Using Lemma 2.3, we can now restrict ourselves by $f\in \Sigma^0$
with
\be\label{2.12}
f(0) \ne b_0.
\end{equation}
Indeed, there exists by this lemma a \qc \ automorphisms $W_\ve$ of
the plane $\hC$, which is conformal outside a domain $E$ located
in a neighborhood of the infinity point, and satisfies
$$
W_\ve(z) = z + a_{2,\ve} z^2 + a_{3, \ve} z^3 + \dots \quad
\text{for} \ \ |z| \le r_1(\ve); \ \ a_{2,\ve} \ne 0, \ \ W_\ve(1) = 1,
W_\ve(\iy) = \iy.
$$
Then
$$
g_\ve(z) = 1/W_\ve(1/z) = z - a_{2,\ve} + b_{1,\ve} z^{-1} + \dots
$$
is conformal outside a small neighborhood of the origin, and
for any $f \in \Sigma^0$ with $f(0) = b_0$ the compozed map
$g_\ve \circ f$ satisfies (2.12).

For $f$ satisfying (2.12), we have by (2.10) $b_1(\mathcal R f) \ne 0$,
and Lemma 2.2 yields
 \be\label{2.13}
\vk(\mathcal R f) = k(\mathcal R f) = k(f), \quad \text{for small} \ \
|t| \le r_{*}(f).
\end{equation}

Put $X_f = (S_f, f(0))$ and consider in $\mathcal F(\T)$ the
homotopy disks
$$
\D(X_f) = \{X_{f_t} = (S_{f_t}, f(0) t): \ |t| < 1\}.
$$
It follows from (2.11) and (2.13) that for all $|t| \le r_{*}$,
$$
c_{\mathcal F(\T)}(\mathbf 0, X_t)
= \tau_{\mathcal F(\T)}(\mathbf 0, X_t)
= d_{\mathcal F(\T)}(\mathbf 0, X_t)
= \tanh^{-1} |t|.
$$
Now, arguing similar to the case $\mathbf{2^0}$, one extends these
equalities to all $|t| < 1$ and then to the initial point
$X = (S_f, f(0)$.

Since the spaces $\T_1$ and $\mathcal F(\T)$ are biholomorphically
equivalent,  the corresponding metrics on $\T_1$ obey similar relations,
which provides that for any point $X \in \T_1$ its distance from
the base point $X_0$ in any invariant (no-expanding) metric is equal
to the \Te \ distance
 \be\label{2.14}
c_{\T_1} (X_0, X) = \tau_{\T_1} (X_0, X) = d_{\T_1} (X_0, X).
\end{equation}
\bk
It remains to establish the equality of distances between two
arbitrary points $X_1, \ X_2$ in $\T_1$.

To reduce this case to (2.14), we uniformize the base point
$X_0 = \D \setminus \{0\}$ (with fixed homotopy class) of this space
by a cyclic parabolic Fuchsian group $\G_0$ acting on the unit disk
(using the universal covering $\pi: \D \to X_0$ with $\pi(0) =
0$) and embed the space $\T_1$ holomorphically into $\T$ via
$$
\T_1 = \T \cap \B(\D^*, \G_0) = \Belt(\D, \G_0)/\sim
$$
(where the equivalence relation commutate with the homotopy of \qc \
homeomorphisms of the surfaces). This preserves all invariant
distances on $\T_1$.

Now, fix a \Be \ coefficient $\mu \in \Belt(\D, \G_0)_1$ so that
$X_1 = w^\nu(X_0)$ as marked surfaces (i.e., with prescribed
homotopy classes) and apply the change rule for \Be \ coefficients:
for any $\mu, \ \nu \in \Belt(\C)_1$, the solutions $w^\mu$ of the
corresponding \Be \ equation
$\partial_{\ov z} w = \mu \partial_z w$ on $\hC$ satisfy
$w^\mu \circ w^\nu = w^{\sigma_\nu(\mu)}$, with
$$
\sigma_\nu(\mu) = (\nu + \mu^*)/(1 + \ov \nu \mu^*),
$$
where
$$
\mu^*(z) = \mu \circ w^\nu(z) \ \ov{\partial_z w^\nu(z)}/\partial_z
w^\nu(z).
$$
Thus, for $\nu$ fixed, $\sigma_\nu(\mu)$ depends holomorphically on
$\mu$ as a map $L_\iy(\C) \to L_\iy(\C)$, which defines a \hol
\ automorphism $\sigma_\mu$ of the ball $\Belt(\D, \G_0)$ preserving
its \Te \ metric. This automorphism is compatible
with \hol \ factorizing projections $\phi_{\T_1}$ and
$\phi_{\T_1^*}$ defining the space $\T_1$ (with base point $X_1$).
Thus $\sigma_\mu$ descends to a \hol \
bijective map $\wh \sigma_\mu$ of the space $\T_1$ onto itself,
which implies the \Te \ isometry
$$
\tau_{\T_1} (\phi_{\T_1}(\mu), \phi_{\T_1}(\nu)) =
\tau_{\T_1}(\phi_{\T_1}(\mathbf 0), \phi_{\T_1}(\sigma_\mu(\nu)),
\quad \nu \in \Belt(f^\mu(\D), f^\mu \G_0 (f^\mu)^{-1})_1,
$$
and similar admissible isometries for the \Ca \ and \Ko \ distances
and their equality.

The case of infinitesimal metrics on $\wt \T$ is investigated in a
similar way, which completes the proof of the theorem.

\bk
\section{Applications to geometric complex analysis}

\subsection{General distortion theorem}
Let $L$ be a bounded oriented quasicircle in the complex plane $\C$
separating the origin and the infinite point, with the interior and
exterior domains $D$ and $D^*$ so that $0 \in D$ and $\iy \in D^*$.
Let $\Sigma(D^*)$ denote the collection of uniavlent functions on
$D^*$ with hydrodynamical normalization at $z = \iy$, and
$\Sigma_k(D^*)$ be its subclass formed by functions $f^\mu$
admitting $k'$-\qc \ extensions to $D$ with $k' \le k$, which we
additionally normalize by $f^\mu(0) = 0$ or $f^\mu(1) = 1$,
and let
$\Sigma^0(D^*) = \bigcup_k \Sigma_k(D^*)$.

Consider on this class a holomorphic (continuous and G\^{a}teaux
$\C$-differentiable) functional $J(f)$, which means that for any $f
\in \Sigma^0(D^*)$ and small $t \in \C$,
$$
J(f + t h) = J(f) + t J_f^\prime(h) + O(t^2), \quad t \to 0,
$$
in the topology of uniform convergence on compact sets in $\D^*$.
Here $J_f^\prime(h)$ is a $\C$-linear functional which is lifted to
the strong (Fr\'{e}chet) derivative of $J$ in the norms  of both
spaces$L_\iy(D)$ and $\B(D^*)$ (cf., e.g., [13, Ch. 3]). Any such
functional $J$ is represented by a complex Borel measure on $\C$ and
extends thereby to all \hol \ functions on $D^*$ (cf. \cite{Sc}).

We shall also use the notation $\wh J(\mu) = J(f^\mu)$; the
functional $\wh J$ is holomorphic on the ball $\Belt(D)_1$.

In the class $\Sigma^0(D^*)$, we have a variational formula
representing the maps with close \Be \ coefficients. This variation
is represented via
$$
\om = H^\mu(z) = z - \fc{1}{\pi} \iint\limits_D \mu(\z) g(\z, z) d
\xi d \eta + O(\|\mu\|_\iy^2) \quad (z \in D^*),
$$
with kernel
$$
g(\z, z) = \fc{1}{\z - z} + g_1(\z),
$$
where $g_1(\z)$ is a rational function determined by the
normalization conditions, and the ratio
\newline $O(\|\mu\|_\iy^2)/\|\mu\|_\iy^2$
remains uniformly bounded on $\hC$ as $\|\mu\|_\iy \to 0$ (see,
e.g., \cite{Kr1}). In particular, in our case,
$$
g(\z, z) = \fc{1}{\z - z} - \fc{1}{1 - z},
$$
and a neighborhood of the identity map $\id$ is filled by the maps
 \be\label{3.1}
\begin{aligned}
f^\mu(z) &= z - \fc{1}{\pi} \iint\limits_\D \mu(\z) \Bigl( \fc{1}{\z
- z} - \fc{1}{\z - 1} \Bigr) d \xi d \eta
+ O(\|\mu\|^2)  \\
&= z - \fc{z - 1}{\pi} \iint\limits_\D \fc{\mu(\z) d \xi d \eta}{(\z
- 1)(\z - z)} + O(\|\mu\|^2) \quad \text{as} \ \ \|\mu\| \to 0.
\end{aligned}
\end{equation}
This yields the corresponding functional derivative of $J$ at the
origin
 \be\label{3.2}
\psi_0(z) = J_{\id}^\prime(g(\id, z)).
\end{equation}
The notation for the kernel $g(\id, z)$ applied here and below
relates to the fact that after the change $\z = f^\nu(Z)$ with $\nu$
running over the set $X$ of Beltrami coefficients from $\Belt(\C)_1$
preserving the domain $D$, one obtains for any fixed finite $z \in
D^*$ a continuous (and real analytic) functional $\check{g}(\nu) =
g(f^\nu, \cdot): \ X \to \C$.

We assume that the functional derivative (3.2) is meromorphic on
$\C$ and has in the domain $D$ only a finite number of the simple
poles (hence $\psi_0$ is integrable over $D$). This holds, in
particular, for the general distortion functionals of the form
 \be\label{3.3}
J(f) := J(f(a); \ f(z_1), f'(z_1), \dots \ , f^{(\a_1)}(z_1); \dots;
f(z_p), f'(z_p), \dots \ , f^{(\a_p)}(z_p))
\end{equation}
with $\grad \wh J(\mathbf 0) \ne 0$, where $z_1, \dots \ , z_p$ are
distinct fixed points in $D^*$ with assigned orders $\a_1, \dots,
\a_p$ and $a$ is a fixed point in the domain of quasiconformality
$D$. In this case, one derives, representing $f$ by (3.1), that
$\psi_0$ is a rational function
$$
\wh J_{\id}^\prime(g(\id, z)) = \fc{\partial \wh J(\mathbf
0)}{\partial \om} g(z, a) + \sum\limits_{j=1}^p
\sum\limits_{s=0}^{\a_j-1} \ \fc{\partial \wh J(\mathbf 0)}{\partial
\om_{j,s}} \fc{d^s}{d\z^s} g(w,\z)|_{w=z,\z=z_s},
$$
where $\om = f(a), \ \om_{j,s} = f^{(s)}(z_j)$.

The above theorems on complex geodesics provide a general distortion
theorem which sheds light on underlying features and, on the other
hand, implies the sharp explicit bounds.

\begin{thm} (i) \ For any functional $J$ of type (3.3) whose range domain
$J(\Sigma^0(D^*))$ has more than two boundary points, there exists a
number $k_0(J) > 0$ such that for all $k \le k_0(J)$, we have
the sharp bound
 \be\label{3.4}
\max\limits_{\|\mu\| \le k} |J(f^\mu) - J(\id)| \le
\max\limits_{|t|= k} |J(f^{t|\psi_0|/\psi_0}) - J(\id)|;
\end{equation}
in other words, the values of $J$ on the ball $\Belt(D)_k = \{\mu
\in \Belt(D)_1: \ \|\mu\|_\iy \le k\}$ are placed in the closed disk
$\D(J(\id), M_k)$ with center at $J(\id)$ and radius $ M_\kp =
\max\limits_{|t|= k} |J(f^{t|\psi_0|/\psi_0}) - J(\id)|$. The
equality occurs only for $\mu = t|\psi_0|/\psi_0$ with $|t| = k$.

(ii) \ Conversely, if a functional $J$ is bounded via (3.4) for $0 <
k \le k_0(J)$ with some $k_0(J) > 0$, then up to rescaling
(multiplying $J$ by a positive constant factor),
 \be\label{3.5}
J(f^\mu) = \mathcal F(S_{f^\mu}) + O(\|\mu\|_\iy^2) \quad \text{as} \ \
\|\mu\|_\iy \to 0,
\end{equation}
where $\mathcal F$ is \hol \ on $\T_1$ and its renormalization
$\wt{\mathcal F}(\vp) = \mathcal F(\vp)/\sup_{\vp\in \T_1}
|\mathcal F(\vp)|$
is the defining map for the disk
$\D(\mu_0)$ as a $c_{\T_1}$-geodesic in the space $\T_1$ with the
base point representing the punctured quasidisk $D \setminus \{a\}$.
\end{thm}

\bk\noindent
{\bf Proof}. For certain specific functionals $J$, such a theorem was
proved in \cite{Kr3}, \cite{Kr5}, \cite{Kr6}.
The proof in the general case follows the same lines and essentially
involves Theorem 1.2.
Without loss of generality, one can assume that $J(\id) = 0$.

First we mention the important facts concerning the projections of norm
$1$ in Banach spaces given in the auxiliary Lemmas 4.1 and 4.2,
the proof of which can be found in \cite{EK}.

Let $V$ be a complex Banach space with norm $\|\cdot\|$ differentiable
on $V \setminus \{0\}$, and suppose that
$$
A(v, w) = \lim\limits_{t \to 0} \fc{\|v + t w\|- \|v\|}{t}
\quad \text{for all} \ \ v \in V \setminus \{0\}, \ w \in V;
$$
for every fixed $v \ne 0$ it is a bounded linear functional on $V$.

\begin{lem} Let $W$ be a non-trivial closed (complex) subspace of $V$,
and let $W^\prime$ be the closed subspace
$$
W^\prime = \{w \in V: \ A(v, w) = 0 \quad \text{for all} \ \
v \in W \setminus \{0\}.
$$
There is a projection $P$ of norm $1$ from $V$ onto $W$ if and
only if $W^\prime$ is a complementary subspace to $W$, that is $W
\oplus W^\prime = V$. Further, if $P$ exists, it is unique and
its kernel is $W^\prime$.
\end{lem}

The next lemma is a straightforward modification of the
corresponding lemmas of Royden \cite{Ro} and of Earle and Kra \cite{EK}.

Let $a_1, a_2, \dots, \ a_n$ be distinguished points of a domain
$D \subseteqq \C$, and let $\vp$ and $\psi$ be $L_1$ functions on
$D$, \hol \ and nonzero for $z \in D \setminus \{a_1, a_2, \dots,
\ a_n\}$. Denote the orders of the functions $\vp$ and $psi$ at the
points $a_j$ by $\a_j$ and $\beta_j$, respectively ($\a_j, \
\beta_j \ge -1$). For real $t$, we consider the function
$$
h(t) = \iint\limits_D |\vp(z) + t \psi(z)| dx dy.
$$

\begin{lem} The function $h(t)$ is differentiable near $t = 0$, and
\be\label{3.6}
h^\prime(0) = \Re \iint\limits_D \psi(z) \fc{\ov{\vp(z)}}{\vp(z)} dx dy.
\end{equation}
Moreover, if
$\a_j \le 2 \beta_j + 1$ for all $j = 1, \dots, \ n$,
then the second derivative $h^{\prime\prime}(0)$ exists.
If $\a_j > 2 \beta_j + 1$ for some $j$, then
\be\label{3.7}
h(t) = h(0) + t h^\prime(0) + \sum\limits_1^n c_j \dl_j(t) +
o(\max_j \dl_j(t)),
\end{equation}
where all the $c_j$ are positive constants, and
\be\label{3.8}
\dl_j = \begin{cases} t^2 \log (1/|t|),  \ \ \ &\a_j = 2 \beta_j + 2,  \\
|t|^{1 + (2 + \beta_j)/(\a_j - \beta_j)},     & \a_j > 2 \beta_j +
2.
\end{cases}
\end{equation}
\end{lem}

Note that the equality (3.6) is obtained by applying Lebesgue's
theorem of dominant convergence, which is possible in view of the
obvious inequality
$$
\Big\vert \fc{|\vp + t \psi| - |\vp|}{t}\Big\vert \le |\psi|.
$$
As for the relations (3.7) and (3.8), for any close subdomain
$E \Subset D \setminus \{a_1, \dots, \ , a_n\}$ the integral
$\iint_E |\vp + t \psi| dx dy$ is an infinitely differentiable function
of $t$, hence the derivation of these relations reduces to estimating
the contributions of integrals of the form
$$
I_j(t) = \iint \Bigl( |z - a_j|^\a + t \om_j(z)| - |z - a_j|^\a - t \Re
\Bigl[\om_j(z)\Bigl(\fc{z - a_j}{|z - a_j}\Bigr)^\a \Bigr] \Bigr) dx dy,
$$
over sufficiently small disks $\{|z - a_j| < r\}$.

\bigskip
Now let $f_0$ be any function in $\Sigma_k(D^*)$ maximizing $|J|$ over
$\Sigma_k(D^*)$ (the existence of such $f_0$ follows from compactness).
We may assume that its \Be \ coefficient $\mu_{f_0}$ is
extremal in its class, i.e.,
$$
\|\mu_{f_0}\|_\iy = \inf \{ \|\mu\|_\iy \le k: \ f^\mu|D^* =
f_0|D^*\}.
$$
Suppose that
\be\label{3.9}
\mu_{f_0} \ne t \mu_0 \quad \text{for some} \ \ t \ \
\text{with} \ \ |t| = k,
\end{equation}
where
\be\label{3.10}
\mu_0 = |\vp_0|/\vp_0.
\end{equation}
We show that for small $k$ this leads to a contradiction.

First we establish the following important property of extremal maps
which sheds light to the main underlying features.
Pick in $A_1(D \setminus \{a\})$ the functions
 \be\label{3.11}
\om_p(z) = z^p - 1 - \psi_0(z), \quad p = 1, 2, \dots \ ;
\quad \rho_a(z) = \fc{a - 1}{(z - 1)(z - a)}.
\end{equation}

\begin{lem} For sufficiently small $k \le k_0(J)$, the extremal \Be \
coefficient $\mu_{f_0}$
is orthogonal to all functions (3.11), i.e.,
 \be\label{3.12}
\langle \mu_{f_0}, \psi_p\rangle_D = 0.
\end{equation}
\end{lem}

\noindent
{\bf Proof}. Note that from (3.1) and (3.2),
 \be\label{3.13}
\wh J(\mu) = - \fc{1}{\pi} \langle \mu, \vp_0\rangle_D +  O(\|\mu\|^2),
\end{equation}
and hence,
\be\label{3.14}
|\wh J(\mu_{f_0})| =
\max_{\|\mu\| \le k} |\wh J(\mu)|
= \fc{k}{\pi} \iint\limits_D |\vp_0| dx dy + O(k^2) =
\|J_{\id}^\prime\| k + O(k^2),
\end{equation}
where $\|J_{\id}^\prime\|$ is defined by (3.1). Consider the
auxiliary functional
 \be\label{3.15}
\wh J_p(\mu) = \wh J(\mu) +
\xi \langle \mu, \psi_p \rangle_D,
\end{equation}
where $p$ is fixed and $\xi \in \C$. Then, similar to (3.14),
 \be\label{3.16}
\max_{\|\mu\| \le k} |\wh J_p(\mu)| = \fc{k}{\pi} \iint\limits_D
|\vp_0(z) + \xi \psi_p(z)| dx dy + O(\|\mu\|^2)
\end{equation}
and the remainder term estimate is independent of $p$.

Using the properties of the norm
$$
h_p(\xi) = \iint\limits_D |\vp_0(z) + \xi \psi_p(z)| dx dy
$$
deduced from Lemmas 3.2 and 3.3 one obtains from (3.14), (3.16)
that for
small $\xi$ there should be $h_p^\prime(0) = 0$, and
 \be\label{3.17}
\max_{\|\mu\| \le k} |\wh J_p(\mu)| = \max_{\|\mu\| \le k}
|\wh J(\mu_{f_0})| + k o_p(\xi) + O_p(k^2 \xi) + O (k^2).
\end{equation}
Now fix $k$ and consider the classes $\Sigma_{\tau k}(D^*), \ \tau < 1$.
Then from (3.17) we have
 \be\label{3.18}
\max_{\|\mu\| \le \tau k } |\wh J_p(\mu)| =
\max_{\|\mu\| \le \tau k} |\wh J(\mu)|
+ \tau o_p(\xi) + O_p(\tau^2 \xi) + O (\tau^2).
\end{equation}
On the other hand, we have from (3.15) that as
$\xi \to 0, \ \tau \to 0$, there should be
$$
\begin{aligned}
|\wh J_p(\tau \mu_{f_0})| &= |\wh J(\tau \mu_{f_0})|
+ \Re \fc{\ov{\wh J(\tau \mu_{f_0})}}{\wh J(\tau \mu_{f_0})}
\tau \xi \langle \mu_{f_0}, \psi_p\rangle + O(\tau^2 \xi^2) \\
&= |\wh J(\tau \mu_0)| + \tau |\xi| \langle \mu_{f_0}, \psi_p\rangle_\D
+ O(\tau^2 \xi^2)
\end{aligned}
$$
with suitable choices of $\xi \to 0$. Comparison with (3.18) implies the
desired equalities (3.12), completing the proof of the lemma.

\bk
We now apply a geodesic \hol \ map $\mathcal F: \T_1 \to \D$ from
Theorem 1.2
defining the disk $\D(\mu_0)$ as $c_{\T_1}$-geodesic; it determines
a hyperbolic isometry between this disk and $\D$. We lift this map
onto $\Belt(\D)_1$ by
$\Ld(\mu) = \mathcal F \circ \phi_{\T_1}(\mu)$
getting
a \hol \ map of this ball onto the disk. The differential of $\Ld$
at $\mu = \mathbf 0$ is a linear operator $P: L_\iy(\D) \to
L_\iy(\D)$ of norm $1$ which is represented in the form
$$
P(\mu) = \beta \langle \mu, \psi_0\rangle_D \ \mu_0.
$$
Let $P(\mu_{f_0}) = \a(\kp) \mu_0$. Since, by assumption, $f_0$ is
not equivalent to $f^{t_0\mu_0}$ with $|t_0| = \kp$, we have
$$
\Big\{\Ld \Bigl(\fc{t}{\kp}\mu_{f_0}\Bigr): \ |t| < 1\Big\}
\subsetneqq \{|t| < 1\}.
$$
Thus, by Schwarz's lemma,
\be\label{3.19}
|a(k)| < k.
\end{equation}
Now consider the function
$$
\nu_0 = \mu_{f_0} - \a(k) \mu_0
$$
which is not equivalent to zero, due to our assumption (4.5).
We show that $\nu_0$ annihilates integrable \hol \ functions on $\D$.

First of all, we have
$$
\langle \nu_0, \vp \rangle_D = 0
$$
for all $\vp$ from the subspace $W^\prime = \langle \psi_p\rangle$
of $A_1(D)$
spanned by (3.11), since $\langle \mu_{f_0}, \vp\rangle_\D = 0$ by
Lemma 3.4 and $\langle \mu_0, \vp\rangle_D = 0$ by Lemma 2.1 applied
to one-dimensional subspace $W = \{\ld \vp_0: \ \ld \in \C\}$.
To establish that
$$
\langle \nu_0, \vp_0 \rangle_\D = 0,
$$
consider the conjugate operator
$$
P^*(\vp) = \langle \mu_0, \vp\rangle_D \vp_0
$$
which maps $L_1(D)$ into $L_1(D)$ and fixes the subspace
$W = (\om_p, \rho_a)$ of $A_1(D \setminus \{a\})$ spanned by
functions (3.11).

The definition of $\nu_0$ implies $P(\nu_0) = 0$, thus
$$
\langle \nu_0, \vp_0\rangle_D = \ld \langle \nu_0, P^* \nu_0\rangle_D
= \langle P \nu_0, \vp_0\rangle_D = 0.
$$
Since the functions $\vp_0, \psi_p, p = 0, 1, \dots$, form a complete
set in the space
$A_1(D \setminus \{a\})$, we have proved that $\nu_0$ is orthogonal to all
$\vp \in A_1(D \setminus \{a\})$, i.e., belongs to the set
$$
A_1(D \setminus \{a\})^\bot = \{\mu \in L_\iy(D): \ \langle \mu, \vp
\rangle_D = 0 \quad \text{for all} \ \ \vp \in A_1(D \setminus \{a\})\}.
$$

Now we use the well-known properties of extremal \qc \ maps
(see e.g., \cite{GL}, \cite{Kr1}, \cite{RS}).
First of all, since $\mu_{f_0}$ is extremal for $f_0$,
$$
\|\mu_{f_0}\|_\iy = \inf \{|\langle \mu_{f_0}, \vp\rangle_D|: \
\vp \in A_1(D\setminus \{a\}), \|\vp\| = 1 \};
$$
moreover, by Hamilton-Krushkal-Reich-Strebel theorem such an
equality is necessary and sufficient for $\mu \in \Belt(D)_1$ to
be extremal for $f^\mu$. Hence, if $f_0$ is extremal in its class,
then for any $\nu \in A_1(D \setminus \{a\})^\bot$,
 \be\label{3.20}
k = \|\mu_{f_0}\|_\iy = \inf \{|\langle \mu_{f_0} + \nu, \vp \rangle_D: \
\vp \in A_1(D \setminus \{a\}), \|\vp\| = 1 \} \le \|\mu_{f_0} + \nu\|_\iy,
\end{equation}
which implies
$$
k \le \|\mu_{f_0} - \nu_0\|_\iy = \|\a(k) \mu_0\|_\iy = \a(k),
$$
in  contradiction to (3.19). This proves the first part of the
theorem.

\bk
The proof of the converse part {\em (ii)}) is much simpler.
Lifting the original functional $J$ with $J(\id) = 0$ to
$$
I(\mu) = \pi^{-1} \circ J(f^\mu): \ \Belt(D)_1 \to \D,
$$
where $\pi$ is a \hol \ universal covering of the domain $V(J) =
J(\Sigma^0(D^*))$ by a disk $\D_a = \{|z| < a\}$ with $\pi(0) = 0,
\pi^\prime(0) = 1$, one obtains
$$
\pi(\z) = \z + O(\z^2), \quad \z \to 0,
$$
with uniform estimate of the remainder for $|\z| < |\z_0|$, which
implies the asymptotic equality (3.5). The covering functional $I$
generates a \hol \ map $\wt I: \ \T_1 \to \D$ so that $I = \wt I \circ
\phi_{\T_1}$, and by (3.4),
$$
\max\limits_{k(f^\mu) \le k} |I(f^\mu)| = k \quad \text{for} \ \
0 < k < k_1(I).
$$
Restricting the covering map $\wt I$ to the extremal disk
$\{\phi_{\T_1}(t \mu_0^*): |t| < 1\} \subset \T_1$ (where $\mu_0^* =
|\psi_0|/\psi_0$), one derives by Schwarz's lemma that
$\wt I(\phi_{\T_1}(t \mu_0^*)) \equiv t$. Thus the
inverse to this map must be $c_\T$-geodesic, which completes the
proof of Theorem 3.1.

\bk
\subsection{A lower estimate for the bound $k_0(J)$} If the
functional $J$ is bounded on the whole class $\Sigma^0(D^*)$, and
$J(\id) = 0, \ \grad J(\id) \ne 0$, one can derive from the proof of
Theorem 3.1 similar to \cite{Kr3} an effective lower bound for
$k_0(J)$:
$$
k \le k_0(J) = \fc{\|J_{\id}^\prime\|}{\|J_{\id}^\prime\| + M(J)
+ 1},
$$
where
$$
\|J_{\id}^\prime\| = \fc{1}{\pi}\|\psi_0\|_1, \quad \quad M(J) =
\sup_{\Sigma^0(D^*)} |J(f)|.
$$

Theorem 3.1 provides new various explicit estimates controlling the
distortion in both conformal and \qc \ domains simultaneously. A
similar theorem is valid also for univalent functions on bounded
quasidisks $D$, for example, for the canonical class $S_k(D)$ of
univalent functions in $D$ normalized by $f(z) = z + c_2 z^2 +
\dots$ near the origin (provided that $z = 0 \in D$) and admitting
$k$-\qc \ extensions to $\hC$ which preserve the infinite point.
Earlier only very special results have been established here (see
\cite{GR}, \cite{Kr1}, \cite{Ku1}, \cite{Ku2}).

\bigskip\noindent
{\bf 3.3. Remark}. It is essential that the distinguished point $a$
is inner, for $a \in \partial D$ the estimate (3.3) can fail and an
additional remainder term  $O(k^2)$ can appear (see, for example,
K\"{u}hnau's description of the domain of values of $f(1)$ on
$\Sigma_k$ in [22, Part 2]).

\section{New phenomena}

\subsection{Rigidity of extremals} The intrinsic connection between the
extremals of the distortion functionals on functions with \qc \
extensions and complex geodesics causes surprising phenomena which
do not appear in the classical theory concerning all univalent
functions. The differences arise from the fact that in problems for
the functions with \qc \ extensions the extremals belong to compact
subsets of $\Sigma^0(D^*)$ (or in other functional classes), while
the maximum on the whole class is attained on the boundary
functions.

\bk The following consequence of Theorem 3.1, part $(ii)$, provides
strong rigidity of extremal maps.

\begin{cor} In any class of univalent functions with $\kp$-\qc \
extension, no function can be simultaneously extremal for different
\hol \ functionals (3.4) unless these functionals have equal
$1$-jets at the origin.
\end{cor}

\subsection{Example: the coefficient problem for functions with \qc \
extensions}
We mention here an improvement in estimating the
Taylor coefficients. Though the Bieberbach conjecture for the
canonical class $S$ of univalent functions  $f(z) = z +
\sum\limits_2^\iy a_n z^n$ in $\D$ has already been proved by de
Brange's theorem, the old coefficient problem remains open for
univalent functions in the disk with \qc \ extensions. The problem
was solved by the author for the functions with sufficiently small
dilatations.

Denote by $S_k(\iy)$ and $S_k(1)$ the classes of $f \in S$
admitting $k$-\qc \ extensions $\wh f$ to $\hC$ normalized by $\wh
f(\iy) = \iy$ and $\wh f(1) = 1$, respectively. Let
 \be\label{4.1}
 f_{1,t}(z) = \fc{z}{(1 - t z)^2}, \quad |z| < 1, \ \ |t| < 1.
\end{equation}
This function can be regarded as a \qc \ counterpart of the
well-known Koebe function which is extremal for many functionals on
$S$ (see, e.g. \cite{Go}, \cite{Po}).

As a special case of Theorem 3.1, we have a complete solution of the
K\"{u}hnau-Niske problem:

\bigskip
\begin{thm} \cite{Kr3} For all $f \in S_k(\iy)$ and all $k \le 1/(n^2 +
1)$,
 \be\label{4.2} |a_n| \le 2 \kp/(n - 1),
\end{equation}
with equality only for the functions
  \be\label{4.3}
f_{n-1,t}(z) = f_{1,t}(z^{n-1})^{1/(n-1)} = z + \fc{2 t}{n - 1} z^n
+ \dots, \quad n = 3, 4, \dots; \ \ |t| = k.
\end{equation}
The estimate (4.2) also holds in the classes $S_k(1)$ with the same
bound for $k$.
\end{thm}

\bk Until now, no estimates have been obtained for arbitrary $k <
1$, unless $n = 2$; in the last case, $|a_2| \le 2 k$ with equality
for the function (4.1) when $|t| = k$.

\bk The rigidity provided by Corollary 4.1 yields that {\em the
function (4.1) cannot maximize $|a_n|$ in $S_k(\iy)$ even for one $k
<1$, unless} $n = 2$. Hence, for all $k < 1$,
 \be\label{4.4}
\max_{f\in S_k(\iy)} |a_n| > n k^{n-1} \quad (n \ge 3).
\end{equation}

\bk For $n = 3$, this estimate was established also by
K\"{u}hnau-Niske \cite{KN} using elliptic integrals.

\bk Comparing the coefficients $a_n$ of $f_{1,t}$ and $f_{n-1,t}$,
one derives from (4.2) and (4.4) the rough bounds for the maximal
value $k_n$ of admissible $\kp$ in (4.2):
$$
\fc{1}{n^2 + 1} \le k_n < \left[ \fc{2}{n(n - 1)}\right]^{1/(n-2)}.
$$

\subsection{Over-normalized functions}
Another remarkable thing in the distortion theory for univalent
functions with \qc \ extension concerns over-determined
normalizations what reveal the intrinsic features of
quasiconformality. The variational problems for such classes are
originated in the 1960s from several points of view. The results
were obtained mainly in terms of inverse extremal functions
$f_0^{-1}$ (see \cite{Kr1}, \cite{BK}, \cite{GoG}, \cite{Re},
\cite{Sh}). But until now no explicit estimates have been
established. We establish here some general explicit bounds for the
functions univalent in the generic disks using complex geodesics in
the space $\T_1$.

\bk Assume that $z = 1$ lies on the common boundary of $D$ and $D^*$
which separates the points $0$ and $\iy$ and denote by
$\Sigma^0(D^*, 1)$ the class of univalent functions in $D^*$ with
\qc \ extensions across $L$ which satisfy
$$
f(z) = z + \const + O(1/z) \ \ \text{near} \ \ z = \iy; \ \ f(1) =
1,
$$
and by $\Sigma_\kappa(D^*, 1)$ its subclasses consisting of
functions with $\kappa$-\qc \ extensions. Fix in $D^*$ a finite
collection of points $e_1, \dots, e_{m-1}, \ m > 1$, and a point
$e_m$ in the complementary domain $D$. Put
$$
e = (e_1, \dots \ , e_m).
$$
We associate with this fixed point set the following subspaces of
$L_1(\C)$: the span $\mathcal L(e)$ of rational functions
 \be\label{4.5}
\rho_s(z) =  \fc{e_s - 1}{(z - 1)(z - e_s)} \quad s = 1, \dots , m,
\end{equation}
the space $A_1(D_e)$ of integrable \hol \ functions in the punctured
domain $D_e = D \setminus \{e_m\}$, and
$$
\mathcal L_0 = \mathcal L(e) \bigoplus \{c \psi_0: \ c \in \C\},
$$
where $\psi_0 = J_{\id}^\prime(g(\id, \cdot))$ for $g(w, \z)$ given
above. Denote by $\Sigma_k(D^*, 1, e)$ the collection of $f \in
\Sigma_k(D^*, 1)$ satisfying
  \be\label{4.6}
f(e_s) = e_s, \ \ s = 1, \dots , m,
\end{equation}
and let $\Sigma^0(D^*, 1, e) = \bigcup_k \Sigma_k(D^*, 1, e)$.

Note that these classes with over-determined normalization contain
nontrivial maps $f^\mu \ne \id$ for any $k < 1$, which is insured by
Lemma 2.3.

\bk
Our aim is to estimate on such over-normalized classes the
functionals
 \be\label{4.7}
J(f) = J(f(z_1), f'(z_1), \dots \ , f^{(\a_1)}(z_1); \dots; f(z_p),
f'(z_p), \dots \ , f^{(\a_p)}(z_p)),
\end{equation}
with $\grad J(\id) \ne 0$, controlling the distortion on the domain
of conformality. To simplify the notations, assume that $ J(\id) =
0$ and $|J(f)| \le 1$ on $\Sigma(D^*)$.

Now one can use only conditional \qc \ variations whose \Be \
coefficients are orthogonal to the rational quadratic differentials
(4.5) corresponding to the fixed points. So the general variational
technique does not work.

From geometric point of view, the equations
 \be\label{4.8}
f(e_1) = e_1, \dots, f(e_m) = e_m
\end{equation}
define an analytic set $\mathcal Z$ of codimension $m$ in the Banach
domain modeling the \uTs \ $\T$. The generic geodesic disks
$$
\{\phi_\T(t \mu): \ \mu \in \Belt(\D)_1, \ t \in \D\}
$$
defining the distortion estimates in Theorem 3.1 can have at most a
finite number common points with this set $\mathcal Z$. This shows
also that the classes $\Sigma_k(D^*, 1, e)$ contains the maps with
dilatations close to $1$.

\bk\subsection {Sharp explicit bounds for small dilatations} The
following distortion bounds involve $L_1$-distance between the
functional derivative $\psi_0$ and span $\mathcal L(e)$; these
bounds hold also in the case when some fixed points $e_s \in
\partial D$.

\bk \begin{thm} For every functional (4.7) and any finite set $e$ of
fixed points defined above, there exists a positive number $k_0(J,
e) < 1$ such that for all $k \le k_0(J, e)$, we have for any
function $f \in \Sigma_k(D^*, 1, e)$  the sharp bound
 \be\label{4.9}
\max\limits_{\|\mu\|\le k} \ |J(f^\mu)| = \max\limits_{|t| = \kappa}
|J(f^{t|\psi_e|/\psi_e})| = d \kappa + O(\kappa^2)
\end{equation}
with uniformly bounded ratio $O(\kappa^2)/\kappa^2$. Here $\kappa =
\kappa(k) < k$ is the bound for \Te \ norms of $f \in \Sigma_k(D^*,
1, e)$ (i.e, regarding these $f$ as the functions of $\Sigma(D^*,
1)$ with omitted restrictions (4.8)),
 \be\label{4.10}
\psi_e = \psi_0 + \sum\limits_1^m \xi_s \rho_s
\end{equation}
with some constants $\xi_1, \dots, \ \xi_m$ and
 \be\label{4.11}
d = \inf_{\psi \in \mathcal L(e)} \ \|\psi_0 - \psi\|_1 =
\end{equation}
These constants are determined (not necessary uniquely) from the
equalities
 \be\label{4.12}
\langle |\psi_e|/\psi_e, \psi \rangle_D = 0 \quad \text{for all} \ \
\psi \in \mathcal L(e); \ \ \langle |\psi_e|/\psi_e, \psi_0
\rangle_D = d
\end{equation}
and from condition that the norm of functional $\langle f_e,
\rho\rangle$ on $L_1$ equals $1$.
\end{thm}

\bk\noindent {\bf Proof}. In view of Theorem 3.1, the functional $J$
can be regarded as a linear functional on the span generated in
$A(D)$ by $\psi_0$. Since the points $e_s$ are preserved, the
variation of $J$ must be orthogonal to $\mathcal L(e)$.

By the Hahn-Banach theorem, there exists a linear functional $l$ on
$L_1(D)$ such that
  \be\label{4.13}
l(\psi) = 0, \ \ \psi \in \mathcal L(e); \ \ l(\psi_0) = d,
\end{equation}
and
$$
\|l\|_{L_1(D)} = \|l\|_{\mathcal L_0} = 1,
$$
and this norm is minimal on the spaces $\mathcal L_0 \subset A_1(D
\setminus \{e_m\}) \subset L_1(D)$. It is generated by $l(\psi) = t
d$ for $\psi \in \mathcal L_0$ and extended to $L_1(D)$ preserving
the norm. Hence, for any other linear functional $\wt l$ on $L_1(D)$
satisfying (4.13) must be $ \|\wt l\|_{\mathcal L_0} \ge 1$. We need
to take $\kappa l(\psi)$.

The functional $l$ is represented on $L_1(D)$ via
$$
l(\psi) = \iint\limits_D \nu_0(z) \psi(z) dx dy , \quad \psi \in
L_1,
$$
with some $\nu_0 \in L_\iy(D)$ so that
  \be\label{4.14}
\iint\limits_D \nu_0(z) \psi(z) dx dy = 0, \ \ \psi \in \mathcal
L(e); \ \ \iint\limits_D \nu_0(z) \psi_0(z) dx dy = d.
\end{equation}
Since the norm of $l$ on the widest space $L_1(D)$ is attained on
its subspace $\mathcal L(e)$, the function $\nu_0$ is of the form
$\nu_0(z) = |\psi_e(z)|/\psi_e(z)$ with integrable \hol \ $\psi_e$
on $D \setminus \{e_m\}$ given by (4.10).

After extending $\nu_0$ by zero to $D^*$, which yields  an extremal
\Be \ coefficient $\nu_0 \in \Belt(D)_1$ for our \hol \ functional
$J$ (with possible pole at the point $e_m$), one can represent the
map $f^{t\nu_0}$ by (4.10) getting from the second equality in
(4.14) and from indicated minimality of $\|l\|$ the estimate (4.12).

However, this $k$-\qc \ map can move the fixed points $e_s$ to
$f_0(e_s) = e_s + O(k^2) = e_s +O(\kappa(k)^2)$, where $f_0 =
f^{\nu_0}$. For correction, one needs to apply additional
$O(k^2)$-\qc \ variation $h_0$ by Lemma 4.3.

It can be shown, using the uniqueness of \Te \ extremal maps
generated by integrable \hol \ quadratic differential that this
$\psi_e$ is unique in $A_1(D \setminus \{e_m\})$.

The assertion on uniform bound for the remainder in (4.9) follows
from the general distortion results for \qc \ maps (see \cite{Kr1}),
completing the proof.

\bk\noindent {\bf Remark 1}. The assumption $f(1) = 1$ can be
replaced by $f(0) = 0 \in D$; then the fixed point $e_m$ must be
chosen to be distinct from the origin.

\subsection{Some applications} In the case when the set $e$
consists of only one point $e_m = \zeta_0$ located in the domain $D$
of quasiconformality (i.e., the functions $f \in \Sigma^0(D^*, 1,
e)$ are normalized hydrodynamically at $\iy$ and $f(1) = 1, \
f(\zeta_0) = \zeta_0$) one can essentially improve Theorem 4.4
getting a global distortion bound of type Theorem 3.1.

Then the equation $f(\zeta_0) = \zeta_0$ defines in the space $\T_1$
an analytic set $\mathcal Z_0$ of codimension $1$, and the generic
geodesic disks joining the origin of $\T_1$ with the points of
$\mathcal Z_0$ are located partially outside of $\mathcal Z_0$. Take
for the points $X = (S_f, z) \in \mathcal Z_0$ the nonsingular \hol
\ disks $\Delta = h(\D)$ in $\T_1$ joining these points with the
origin and pull back the hyperbolic metric $\lambda_\D$ onto such
disks using the functions $\wh J \circ h: \D \to \D$, where $\wh J$
denotes the lifting of the functional $J$ onto the space $\T_1$.
This yields a \sh \ metric $\ld_J$ on $\Delta$ of curvature $-4$ at
its noncritical points. In particular, there is a domain $\Omega
\subset \Delta$ containing the origin filled by the points of
$\mathcal Z_0$. Comparing this metric with the infinitesimal Finsler
metric $F_\T(S_f, v)$ (which defines the \Te \ dilatation $\kappa =
\tanh \tau_\T$) similar to part $\mathbf 3^0$ of the proof of
Theorem 3.1, one derives the following bound:

\begin{prop} For $0< k < k_0(J, \zeta_0) \le 1$,
 \be\label{4.15}
\max \{|J(f^\mu): \ f^\mu \in \Sigma(D^*, 1, e), \ \|\mu\|\le k \} <
d \kappa,
\end{equation}
where $\kappa = \kappa(k)$ and $d$ is the $L_1$-distance between the
differential $\psi_0$ and the line
$$
c \fc{\zeta_0 - 1}{(z - 1)(z - \zeta_0)}, \quad c \in \C.
$$
\end{prop}

In contrast to Theorems 3.1 and 4.3, the bound (4.15) is not sharp.

\bk Similar results are valid also for the over-normalized functions
in bounded quasidisks. We illustrate those again on the coefficient
problem:

{\em Find $\max |a_n| \ (n \ge 2)$ for the functions $f \in
S_k(\iy)$ leaving a given set $e = (e_1, \dots, \ e_m) \subset
\partial (\D^* \setminus \{\iy\})$ fixed} (hence some $e_s = \iy$).

\bk
\begin{thm}
For any $n \ge 2$, there is a number $k_n(e) < 1$ such that for $k
\le k_n$ and all $f \in S_k(\iy)$, which fix a given set $e = (e_1,
\dots, \ e_m)$, we have the sharp bound
 \be\label{4.16}
\max\limits_{\|\mu\|\le k}  |a_n(f^\mu)| = \max\limits_{|t| =
\kappa} |a_n(f^{t|\psi_n|/\psi_n})| = \fc{2d_n \kappa}{n-1} +
O(\kappa^2),
\end{equation}
where $\kappa(k)$ is again the bound for \Te \ norms, and similar to
(4.10),(4.11),
$$
\psi_n(z) = c z^{-n-1} + \sum\limits_1^m \xi_s \rho_s(z), \quad d_n
= \inf_{\mathcal L(e)} \|\psi_n - \psi\|_1.
$$
The constants $c, \xi_s$ are determined from the equations of type
(4.12), and the remainder in (4.16) is estimated uniformly for all
$k \le k_n$.
\end{thm}

\bk\noindent {\bf Remark 2}. The distortion bounds of type (4.9)
given by Theorem 4.4 and its corollaries hold in somewhat weakened
form (up to terms $O(k^2)$) for the maps preserving an infinite
subset $e$ in $D$, provided that the corresponding class
$\Sigma_k(D^*, 1, e)$ contains the functions $f^\mu \ne \id$.

The proof is similar but now the quadratic \hol \ differentials
$\psi_e$ defining the extremal functions are represented instead of
(4.10) in the form
$$
\psi_e = \psi_0 + \psi, \quad \psi \in \mathcal L(e)
$$
and there are no variations of type Lemma 4.3 for the infinite sets.

\subsection{Shift-like theorems}
The above theorems also make a progress in solving old question
related to Teichm\"{u}ller's Verschiebungssatz of 1944. This theorem
provides explicitly the extremal map among \qc \ automorphisms of
the disk $\D$ preserving fixed all boundary points and moving a
given point $z_0 \in \D$ into a given $w_0$. It found interesting
applications. There were given some its modifications, but until now
no finite dimensional analogs (i.e., for finite sets of fixed
boundary points) have been established. Theorems 3.1 and 4.4 give a
particular answer to this question.

Consider \hol \ functionals $J(f)$ of the form
 \be\label{4.17}
J(f) = J(f(z_1), \dots, f(z_n))
\end{equation}
depending on the values of maps $f \in \Sigma(D^*)$ in the
distinguished points $z_j \in D \cup D^*$. Let again $J(\id) = 0, \
\grad J(\id) \ne 0$. For any fixed finite set $e = (e_1, \dots, e_m)
\subset \partial \D^*$, the maximal value of $|J(f)|$ on the class
$\Sigma(\D^*, e)$ is obtained via (4.9) and growths monotonically
with the dilatation $k$. Thus, one can prescribe $r <
\max_{\Sigma(\D^*, e)} |J(f)|$ and ask on the minimal dilatation $k$
on which such $r$ is attained on the maps $f^\mu \in \Sigma(D^*)$
preserving the points of the set $e$. Theorem 4.4 yields the
following geometric result.

\begin{thm} For small $r > 0$, the minimal dilatation, on which
any functional (4.17) attains its level surface $L_r = \{|J(f)| =
r\}$, equals
$$
\kappa = r/d + O(r^2)
$$
with $d$ defined by (4.11). This estimate is sharp.
\end{thm}

The problem of establishing sharp explicit global distortion
estimates for over-normalized maps, even such as Theorem 3.1, is
very complicated and remains open.

\bk I am thankful to the referee for his comments and suggestions.

\bk
\medskip
{\small\em{ \leftline{Department of Mathematics, Bar-Ilan
University} \leftline{5290002 Ramat-Gan, Israel}
\leftline{and Department of Mathematics, University of Virginia,}
\leftline{Charlottesville, VA 22904-4137, USA}}}

\end{document}